\documentclass[11pt, a4paper]{article}
\usepackage[utf8]{inputenc} 
\usepackage[frenchb, german, english]{babel}
\usepackage[T1]{fontenc}		
\usepackage{url}
\usepackage{amssymb}		
\usepackage{epsf}			
\usepackage{amsmath}		
\usepackage{newtxmath,newtxtext}
\usepackage{graphicx} 		
\usepackage[bottom]{footmisc}	

\title{Diophantus redivivus:\\ is Diophantus an early-modern classic?}
\author{Catherine Goldstein,\\ 
CNRS, Institut de mathématiques de Jussieu-Paris Rive Gauche,\\ UMR 7586, Sorbonne Université, Université de Paris,\\  Case 247, 4 place Jussieu, 75252 Paris Cedex 05, France,\\ catherine.goldstein@imj-prg.fr}
\date{June 17, 2021}
\begin{document}

 \maketitle

Many early-modern mathematical books incorporated at least a part of Diophantus' \emph{Arithmetica},  from Jacques de Billy's \emph{Diophanti Redivi Pars prior et posterior} to John Kersey's \emph{Third and Fourth Books of the Elements of algebra}  or Jacques Ozanam's \emph{Recr\'eations math\'ematiques}. Diophantine questions regularly circulated among mathematicians of the time in the context of exchanges of problems or challenges \cite{CGGoldstein2013}.  It is thus tempting to consider Diophantus's opus magnum as a classic. However, I argued in my talk that, while Diophantus was indeed a classical author for early-modern mathematicians, his main work did not become a classical book.

The first point is easy to establish. According to the definition of the \emph{Dictionnaire de l'Acad\'emie} \cite[vol.~I, p.~197]{CGDictionnaire1694}, 
\begin{quote}
``classical'' [is] only in use in this sentence: \emph{a classical author}. That is: an ancient and much approved Author who is an authority on the subject he deals with. \emph{Aristotle, Plato, Livy, \&c. are classical Authors}.
\end{quote}
Statements  establishing such a status for Diophantus abound. To give just one example among many,  in his 1660 \emph{Diophantus geometra}, Billy claims \cite[Lectori Benevolo\dots]{CGBilly1660}: 
\begin{quote}
Who does not know Diophantus does not deserve the name of mathematician ; indeed, what Cicero is among the Orators, Virgil among the Poets, Aristotle among the Philosophers, Saint Thomas among the Theologists, Hippocrates among the Physicians, Justinian among the Lawmakers, Ptolemy among the Astronomers, Euclid among the Geometers, is the very same as what Diophantus is among the Arithmeticians ; he who surpasses all the others by a long interval is their coryphaeus, and easily their prince.
\end{quote}

Still, several historiographical issues are at stake. First of all, at least two authoritative versions of Diophantus' \textit{Arithmetica} are referred to by most mathematicians, 
Franciscus Vieta's \emph{Zetetica} \cite{CGViete1591} and the heavily commented edition, with a Latin translation, published by Claude-Gaspard Bachet de Meziriac \cite{CGBachet1621}. 
These texts belong to different genres and provide their readers with different organizations and selections of Diophantine material, different symbolisms and textual marks, and even different  ideas of what constitutes an adequate solution \cite{CGMorse1982} (see also Abram Kaplan's contribution to this workshop). 

On the other hand, both, in various degrees, present Diophantus as the father of algebra, a clue followed by  most early-modern mathematicians who followed them. Billy, mentioned above, for instance, goes on  \cite[Lectori Benevolo\dots]{CGBilly1660}: 
\begin{quote}
Diophantus remained within the limits of Arithmetic, not the vulgar kind that is taught to children and merchants, but another, more subtle one, that one calls Algebra and that is the science of unknown numbers starting from hypotheses.
\end{quote}  

Even this agreement poses problems: one of its consequences is that Diophantine questions are usually treated in textbooks on algebra, another that algebraic  formulas are often presented as a natural and valid generalization of Diophantus's original, unique, solution in fractions. Some, like Fermat, advocated at the time against this trend, but in favor of a search for integer solutions only. However, none of these developments correspond to the current idea about Diophantine questions,  with its emphasis on the description and computation of rational solutions \cite{CGGoldstein1995, CGSchappacher1998}. 

A last,   intriguing, issue concerns the French scene particularly. In French literature, ``classical'' has been more and more associated with the new standards emphasized in the framework of, or parallel to, the courtly culture of Louis XIV's times; opposed to the heavy and contrived volumes of the schools \cite{CGDhombres1996}, classical texts, in this sense, were supposed to share a set of characteristics, such as elegance, simplicity, naturalness and correctness \cite{CGViala1993, CGRibot2018}. A question is thus  whether  Diophantine problems found a home in mathematical books that can be qualified as classical, and how. 

To study these issues, I analyzed in the talk several works linked to Diophantus's \emph{Arithmetica}, focussing for the sake of time to  French authors of the second half of the seventeenth century: Billy's \emph{Diophantus geometra} \cite{CGBilly1660} and \emph{Diophantus redivivus} \cite{CGBilly1670, CGCassinet1987},  Ozanam's manuscript \emph{Diophante reduit \`a la specieuse} and \emph{Trait\'e des simples, des doubles and des triples \'egalit\'es pour la solution des probl\`emes en nombres} \cite{CGCassinet1986}, as well as his  \emph{Nouveaux elemens d'algebre} \cite{CGOzanam1702}, Jean Prestet's \emph{Elemens de mathematiques} (both editions) \cite{CGPrestet1675, CGPrestet1689, CGAsselah2005}, and finally Bernard Frenicle de Bessy's posthumous \emph{Trait\'e des triangles rectangles en nombres} \cite{CGGoldstein1995}. For each of them, I surveyed  the selection and organization of the material, as well as the individual treatment of some of the problems and their solutions. This shows that if Diophantus' \emph{Arithmetica} remained an inspiration for  new isolated mathematical problems, in algebra, number theory and even in Euclidean geometry, they entered into a large variety of partial reconstructions, within different genres. There was no agreement at the time on what was a satisfactory solution to Diophantine problems: even with algebraic formulas, for instance, and the claims of their authors to go further than Diophantus in providing infinitely many solutions instead of a single one, no proof of the infinity of the solutions was ever given, nor that they were all obtained. Attempts to restructure Diophantine problems and eventually integrate them in a text adapted to the new fashion (such as Ozanam's \emph{Trait\'e des simples, des doubles and des triples \'egalit\'es}  or Frenicle's  \emph{Trait\'e des triangles rectangles en nombres}) either were not published in their original form or relied on a severe selection of a few problems and topics. 

Contrary to other authors, thus, there was no Diophantus for the \textit{honn\^ete homme}. And early-modern Diophantus  appears as classical author without a classical text.

\end{document}